\documentstyle[12pt]{article}
\def\Dslash{\not\hspace{-0.7ex}D_h\hspace{0.7ex}}

\def\bbbr{{\rm I\!R}} 

\def\bbbc{{\mathchoice {\setbox0=\hbox{$\displaystyle\rm C$}\hbox{\hbox
to0pt{\kern0.4\wd0\vrule height0.9\ht0\hss}\box0}}
{\setbox0=\hbox{$\textstyle\rm C$}\hbox{\hbox
to0pt{\kern0.4\wd0\vrule height0.9\ht0\hss}\box0}}
{\setbox0=\hbox{$\scriptstyle\rm C$}\hbox{\hbox
to0pt{\kern0.4\wd0\vrule height0.9\ht0\hss}\box0}}
{\setbox0=\hbox{$\scriptscriptstyle\rm C$}\hbox{\hbox
to0pt{\kern0.4\wd0\vrule height0.9\ht0\hss}\box0}}}}
\def\bbbq{{\mathchoice {\setbox0=\hbox{$\displaystyle\rm Q$}\hbox{\raise
0.15\ht0\hbox to0pt{\kern0.4\wd0\vrule height0.8\ht0\hss}\box0}}
{\setbox0=\hbox{$\textstyle\rm Q$}\hbox{\raise
0.15\ht0\hbox to0pt{\kern0.4\wd0\vrule height0.8\ht0\hss}\box0}}
{\setbox0=\hbox{$\scriptstyle\rm Q$}\hbox{\raise
0.15\ht0\hbox to0pt{\kern0.4\wd0\vrule height0.7\ht0\hss}\box0}}
{\setbox0=\hbox{$\scriptscriptstyle\rm Q$}\hbox{\raise
0.15\ht0\hbox to0pt{\kern0.4\wd0\vrule height0.7\ht0\hss}\box0}}}}
\def\bbbt{{\mathchoice {\setbox0=\hbox{$\displaystyle\rm
T$}\hbox{\hbox to0pt{\kern0.3\wd0\vrule height0.9\ht0\hss}\box0}}
{\setbox0=\hbox{$\textstyle\rm T$}\hbox{\hbox
to0pt{\kern0.3\wd0\vrule height0.9\ht0\hss}\box0}}
{\setbox0=\hbox{$\scriptstyle\rm T$}\hbox{\hbox
to0pt{\kern0.3\wd0\vrule height0.9\ht0\hss}\box0}}
{\setbox0=\hbox{$\scriptscriptstyle\rm T$}\hbox{\hbox
to0pt{\kern0.3\wd0\vrule height0.9\ht0\hss}\box0}}}}
\def\bbbs{{\mathchoice
{\setbox0=\hbox{$\displaystyle     \rm S$}\hbox{\raise0.5\ht0\hbox
to0pt{\kern0.35\wd0\vrule height0.45\ht0\hss}\hbox
to0pt{\kern0.55\wd0\vrule height0.5\ht0\hss}\box0}}
{\setbox0=\hbox{$\textstyle        \rm S$}\hbox{\raise0.5\ht0\hbox
to0pt{\kern0.35\wd0\vrule height0.45\ht0\hss}\hbox
to0pt{\kern0.55\wd0\vrule height0.5\ht0\hss}\box0}}
{\setbox0=\hbox{$\scriptstyle      \rm S$}\hbox{\raise0.5\ht0\hbox
to0pt{\kern0.35\wd0\vrule height0.45\ht0\hss}\raise0.05\ht0\hbox
to0pt{\kern0.5\wd0\vrule height0.45\ht0\hss}\box0}}
{\setbox0=\hbox{$\scriptscriptstyle\rm S$}\hbox{\raise0.5\ht0\hbox
to0pt{\kern0.4\wd0\vrule height0.45\ht0\hss}\raise0.05\ht0\hbox
to0pt{\kern0.55\wd0\vrule height0.45\ht0\hss}\box0}}}}

\begin{document}

\rightline{DAMTP/97-7}

\begin{center}
{\Large {\bf Twisted $h$-spacetimes and invariant equations}}
\end{center}


\begin{center}
{\large{J.A. de Azc\'{a}rraga$^{\dagger}$\footnote{St. John's College Overseas
Visiting Scholar. On sabbatical leave of absence from the 
Department of F\'{\i}sica Te\'{o}rica and IFIC,
Centro Mixto Universidad de Valencia-CSIC,
E-46100-Burjassot (Valencia), Spain.}, 
P.P. Kulish$^{\ddagger}$ and F. Rodenas$^{\ast}$}}
\end{center}

\noindent
{\small {\it $ \dagger$ Department of Applied Mathematics and Theoretical
Physics, Silver St., Cambridge CB3 9EW, U.K.}}


\noindent
{\small {\it $ \ddagger$ Steklov Mathematical Institute, Fontanka 27, 
191011-St. Petersburg,  Russia.}} 

\noindent
{\small{\it $ \ast $ Departamento de Matem\'atica Aplicada,
ETSA,  Universidad Polit\'ecnica de Valencia
E-46071 Valencia, Spain.}}

\vspace{1\baselineskip}

\begin{abstract}

    We analyze the  $h$-deformations of the Lorentz group and their 
associated spacetimes. We prove that they have a twisted character and 
give explicitly the twisting matrices. After studying the representations
of one of the deformed spacetime algebras, we discuss the Klein-Gordon
operator. It is found 
that the $h$-deformed d'Alembertian has plane wave solutions of the same
form as the standard ones. We also give explicit expressions for the 
$h$-gamma matrices defining the associated  Dirac equations.   
\end{abstract}

\section{Introduction}

\indent

It is about fifty years since the first attempts to discretize 
the spacetime manifold were probably first made in \cite{SNYDER}, where
the introduction of a smallest unit of length in spacetime led to 
non-commutative spacetime coordinates. The quantum structure of spacetime is 
still  being  discussed from different viewpoints \cite{DOPLI} (see also
\cite{KEMPF} and references therein) which quite often turn out to be 
associated with various aspects  of  non-commutative geometry  \cite{CONNES}.
This has been used very recently to propose an action unifying gravity 
and the standard model at very high energies \cite{CC}.
  A direct introduction of a spacetime lattice structure breaks 
the relativistic invariance of the theory, but  
the recent generalization of the Lie groups to deformed (`quantum') 
groups gives rise to  
another possibility of constructing  non-commutative spacetime 
coordinates in order to find  the algebraic analogues
of the special relativity relations.

The initial step of the quantum group approach consists in selecting 
a deformed Lorentz group and then defining  a quantum Minkowski spacetime  
algebra by using covariance arguments similar to those used to define the 
`quantum' plane \cite{MAN,WessZ}. Although the best known and 
well developed deformation with dimensionless parameter 
\cite{WA-ZPC48,OSWZ-CMP}  
is based on  the $q$-deformed $SL_q(2)$
group,  there are a variety of other deformations  and 
even a classification scheme \cite{WOZA,POWO,JPA96} 
(see also \cite{LCAST}). However in the papers mentioned above 
different techniques are used which do not
 permit an easy comparison of the results obtained.

The aim of this paper is to provide a detailed analysis of 
$h$-deformed Minkowski spacetimes  and their properties, which are  
related to the  twisted $SL_h(2)$ \cite{JPA96,BANACH}.
The $h$-deformation is specially simple due to the natural $*$-structure
\cite{PODLES,BANACH} of the $h$-spacetime algebras, which is absent in general.
Recently, and in the framework of \cite{WOZA,POWO}, a formal solution 
to the deformed Klein-Gordon and Dirac equation has been proposed 
\cite{PODLES} by 
introducing `non-commutative eigenvalues' of the deformed spacetime
derivatives. However, the introduction of  simplifying commutation 
relations for deformed spacetime coordinates and derivatives 
was done in a non-covariant way: the set of commutation 
relations of coordinates and derivatives with `non-commutative 
eigenvalues' is not preserved under deformed Lorentz group 
transformation (coaction). We shall consider also this problem in our
framework and find a solution consistent with covariance. The problem 
of the realization of Dirac matrices will also be  solved naturally 
for the twisted deformation. 

The paper is organized as follows.
In Sec. 2 we recall the properties of the $h$-deformed Minkowski 
spacetime \cite{JPA96,BANACH} following the $R$-matrix 
formalism \cite{AKR,FTUV94-21}. The different deformed Lorentz 
groups and corresponding Minkowski spaces as well as their properties 
are directly related to a set of four $R$-matrices $R^{(j)}$ 
($j$=1,2,3,4), solutions to the Yang-Baxter equation (YBE) and FRT-relations.  
Selecting the $R$-matrices associated with $SL_h(2)$ \cite{JPA96}
we are able to determine a  complete and unique set 
of covariant commutation relations among coordinates, derivatives and 
one-forms.

Sec. 3 gives the explicit expressions for the metric matrices for 
coordinates and derivatives, and shows how the mild (twisted) nature of 
the deformation is also encoded in the dilatation operator. Sec. 4 describes
the $h$-deformed spacetime algebras and  representations. Sec. 5
includes a general discussion of the Lorentz group deformations connected 
with triangular $R$-matrices and shows how the twisted nature of $SL_h(2)$
is also incorporated in the $h$-Lorentz groups.
The properties of the resulting ${\cal R}$-matrix (16$\times$16) are derived 
as consequences of the corresponding properties of its components,
triangularity and twisting character, and the twisting  matrices are
constructed.
Using the exchange ${\cal R}$-matrix and a general approach to the 
construction of deformed 
Dirac operators \cite{FTUV94-21}, the problem of finding the
$h$-Dirac matrices  is solved in Sec. 6. 
Sec. 7 is devoted to the study of the solutions of the $h$-deformed 
Klein-Gordon equation in a way consistent with covariance. Finally, 
the paper closes with some conclusions.

\setcounter{equation}{0}

\section{$h$-deformed Minkowski algebras and differential operators}

\indent

We review first the properties of the two
deformed   Minkowski spaces \cite{JPA96,BANACH}  associated with the 
`Jordanian' or `non-standard' $h$-deformation, $SL_h(2)$, of $SL(2,\bbbc)$
\cite{DEMI,EOW,ZAKRZ,KAR}. The $GL_h(2)$  deformation is 
defined as the associative algebra generated by the 
entries  $a,b,c,d$ of a matrix $M$, the commutation properties of which
may be  expressed by an `FRT'  equation,
$R_{12}M_1M_2=M_2M_1R_{12}$  \cite{FRT1}, where  $R$  is
the triangular solution of the Yang-Baxter equation 
\begin{equation}\label{1.5}
R_h= \left[ 
\begin{array}{cccc}
1 & -h & h & h^2 \\
0 & 1 & 0 & -h \\
0 & 0 & 1 & h \\
0 & 0 & 0 & 1 
\end{array} \right] \;,
\; \hat{R}_h \equiv {\cal P}R_h= \left[ 
\begin{array}{cccc}
1 & -h & h & h^2 \\
0 & 0 & 1 & h \\
0 & 1 & 0 & -h \\
0 & 0 & 0 & 1 
\end{array} \right] \;,   \; {\cal P}R_h{\cal P}=R_h^{-1} \;,
\end{equation}	
\noindent
and  ${\cal P}$ is the permutation operator (${\cal P}_{ij,kl}=
\delta_{il} \delta_{jk}$, ${\cal P}={\cal P}^{-1}$,
${\cal P}(A \otimes B) {\cal P}= B \otimes A$ if the entries of $A$ and $B$
commute).
The commutation relations of the  algebra  generators in $M$ are 
\begin{equation}\label{1.6}
\begin{array}{lll}
{} [a,b]= h(\xi-a^2) \quad, & \quad [a,c]=hc^2 \quad,&\quad 
[a,d]= hc(d-a) \quad, \\
{}   [b,c]=h(ac+cd) \quad,  & \quad
 [b,d]=h(d^2-\xi) \quad, & \quad [c,d]=-hc^2 \quad ;
\end{array}
\end{equation}	
\begin{equation}\label{1.7}
\xi \equiv det_hM = ad - cb - h cd  \quad.
\end{equation}	
\noindent
Setting $\xi =1$  reduces  $GL_h(2)$ to $SL_h(2)$.
The matrix $\hat{R}_h$ (being triangular, $\hat{R}_h^2=I_4$) has two 
eigenvalues ($1$ and $-1$) and a spectral decomposition in terms 
of a rank three projector $P_{h+}$ and a rank one projector $P_{h-}$,
\begin{equation}\label{1.7.2}
\hat{R}_h  =P_{h+} - P_{h-} \quad, \quad  P_{h \pm} \hat{R}_h = \pm P_{h \pm} 
\quad, \quad  P_{h \pm}= \frac{1}{2}(I \pm \hat{R}_h)         \quad. 
\end{equation}
\noindent
 The deformed determinant $\xi$
in (\ref{1.7}) may be then expressed as
$(det_h M)P_{h-} := P_{h-} M_{1} M_{2}$. 
The following relations  for the $h$-symplectic metric
have an obvious  equivalent in the undeformed  case:
\begin{equation}\label{1.7.4}
\epsilon_hM^t\epsilon_h^{-1}=M^{-1} \,, \, \epsilon_h= \left(
\begin{array}{cc}
      h & 1 \\
      -1 & 0
      \end{array} \right) \,, \,  \epsilon_h^{-1}= \left(
\begin{array}{cc}
      0 & -1 \\
      1 & h
      \end{array} \right) \,, \, 
P_{h- \;ij,kl}=\frac{-1}{2}\epsilon_{h\;ij}\epsilon^{-1}_{h\;kl} \,.
\end{equation}
\noindent
Deformed `groups' related with  
different values of $h \in \bbbc $
are equivalent and  their $R_h$ matrices  are related by  a similarity 
transformation.  Thus, from now on, we shall take  
$h \in \bbbr$.

 The determination of a complete set of 
deformations of  the Lorentz group \cite{WOZA,JPA96} requires replacing 
 \cite{POWO,WA-ZPC48,OSWZ-CMP}
the $SL(2,\bbbc)$ matrices $A$ in the undeformed expression $K'=AKA^{\dagger}$ 
($K=K^{\dagger}=\sigma_{\mu}x^{\mu}$) by the generator 
matrix $M$ of a deformation of $SL(2,\bbbc)$,  and 
the characterization of all possible commutation relations
among the generators 
($a,b,c,d$) of $M$ and ($a^*,b^*,c^*,d^*$) of $M^{\dagger}$.
In particular, for  the deformed Lorentz  groups associated with $SL_h(2)$,
the $R$-matrix form of these commutation relations  may be   expressed  by
\begin{equation}\label{2.3}
\begin{array}{ll}
R_h M_1M_2=M_2M_1 R_h \quad, & \quad
M_1^{\dagger}R^{(2)}M_2=M_2R^{(2)}M_1^{\dagger} \quad, \\
M_2^{\dagger}R^{(3)}M_1=M_1R^{(3)}M_2^{\dagger} \quad, & \quad
R^{\dagger}_h M_1^{\dagger}M_2^{\dagger}
=M_2^{\dagger}M_1^{\dagger}R^{\dagger}_h \quad,
\end{array}
\end{equation}
\noindent
where $R^{(3) \,\dagger}=R^{(2)}={\cal P}R^{(3)}{\cal P}$
(`reality' condition for $R^{(3)}$). 
The consistency of these relations is assured if 
$R^{(3)}$ satisfies the FRT
equation  (see \cite{FTUV94-21,FM} in this respect)
\begin{equation}\label{2.4}
R_{h\, 12}R^{(3)}_{13}R^{(3)}_{23}=R^{(3)}_{23}R^{(3)}_{13}R_{h\,12}
\quad.
\end{equation}
\noindent
This equation, considered as an FRT equation, 
indicates that $R^{(3)}$ is a representation of  $GL_h(2)$,  
$(M_{ij})_{\alpha \beta}=
R^{(3)}_{i \alpha , j \beta}$. 
The solutions of these equations\footnote{ The first
sentence in \cite{JPA96} in the Remark after (36) should be deleted.}
\cite{JPA96,WOZA} ($h,r \in \bbbr$),
\begin{equation}\label{4.8}
\mbox{{\bf 1.}} \quad  R^{(3)}=\left[
\begin{array}{cccc}
1 & 0 & 0 & 0  \\
0  & 1 & r  & 0   \\
0  & 0 & 1 & 0 \\
0  & 0 & 0 & 1
\end{array} \right] \quad ,   \quad
\mbox{{\bf 2.}} \quad   
R^{(3)}=\left[
\begin{array}{cccc}
1 & 0 & -h & 0  \\
-h  & 1 & 0  & h   \\
0  & 0 & 1 & 0 \\
0  & 0 & h & 1
\end{array} \right] \;, 
\end{equation}
\noindent
characterize the two  $h$-deformed Lorentz groups, which
will be denoted $L_h^{(1)}$ and $L_h^{(2)}$ respectively.

\vspace{1\baselineskip}

To introduce the deformed Minkowski {\it algebra} ${\cal M}^{(j)}_h$ 
associated with the $h$-Lorentz group  $L^{(j)}_h$  ($j$=1,2) we extend   
$K'=AKA^{\dagger}$ above
by stating that in  the deformed case   the corresponding 
$K^{(j)}$ generates  a comodule algebra for the coaction $\phi$ defined by
\begin{equation}\label{2.5}
\phi : K \longmapsto
K' = M KM^{\dagger} \;, \; K'_{is} = M_{ij} M^{\dagger}_{ls}
K_{jl} \;, \; K= \left( 
\begin{array}{cc}
\alpha & \beta \\
\gamma & \delta  \end{array} \right) 
= K^{\dagger} \;,\; \Lambda = M \otimes M^* \;,
\end{equation}
\noindent
where $\Lambda$ is the $L_h^{(j)}$ matrix and, as usual,  it is assumed 
that the elements of $K$, 
which now do not commute among themselves, commute with those of $M$
and $M^{\dagger}$.  
We now demand that the commuting 
properties of the entries  of $K$  are preserved by
(\ref{2.5}).  The use of covariance arguments 
to characterize the algebra generated by the elements of  $K$
has been extensively used, and the resulting equations are associated 
with the name of reflection equations (see \cite{K-SKL,KS} and references 
therein) or, in a more
general setting, braided algebras \cite{MAJ-LNM,SMBR2}. In  the present
$SL_h(2)$ case, the 
commutation properties of the entries of the hermitian  matrix $K$ generating 
a deformed Minkowski algebra ${\cal M}_h$ are given by 
\begin{equation}\label{2.6}
R_h K_{1} R^{(2)} K_{2} = K_{2} R^{(3)} K_{1} R^{\dagger}_h \quad,
\end{equation}
\noindent
where the $R^{(3)}=R^{(2)\, \dagger}$ matrices are those given in 
(\ref{4.8}); using (\ref{2.3}), it is easy to check that (\ref{2.6}) is
preserved under $\phi$ in (\ref{2.5}).

The deformed  Minkowski length and metric, invariant under 
a Lorentz transformation
(\ref{2.5}) of $L^{(j)}_h$, is defined \cite{JPA96} through the $h$-deformed
determinant of $K$
\begin{equation}\label{gdet}
(det_{h}K)P_{h-}P_{h-}^{\dagger} 
= - P_{h-} K_1 \hat{R}^{(3)}K_1 P_{h-}^{\dagger} \quad,
\end{equation}
\noindent
where $P_{h-}P_{h-}^{\dagger}$ is a projector since
$(P_{h-}P_{h-}^{\dagger})^2=
\left( \frac{2+h^2}{2} \right)^2 P_{h-}P_{h-}^{\dagger}$.
$det_{h}K$ is invariant, central and, since $\hat{R}^{(3)}$ 
and $K$ are hermitian, real; thus, it defines  the 
{\it deformed Minkowski length} $l_h$
for the $h$-deformed spacetimes ${\cal M}^{(j)}_h$.

To describe the differential calculus on $h$-Minkowski spaces,
we need expressing the different commutation relations among the fundamental 
objects: deformed coordinates, derivatives and  one-forms.
Following the approach of \cite{AKR,FTUV94-21,KARPACZ} 
to the differential calculus on 
Minkowski algebras associated with the standard deformation $SL_q(2)$, we 
introduce the reflection equations expressing the  commutation
relations defining the algebras of $h$-derivatives and $h$-one-forms
($h$-differential calculi have been considered in \cite{KAR} and 
in \cite{PODLES}  for quantum $N$-dimensional homogeneous spaces).

The  derivatives are expressed in terms  of an object 
$Y$ transforming  covariantly (cf. (\ref{2.5})) {\it i.e.},

\begin{equation}\label{dif1}
Y \longmapsto Y'=(M^{\dagger})^{-1} Y M^{-1} \quad , \quad
Y= \left[ \begin{array}{cc}
                  \partial_\alpha & \partial_\gamma  \\
                 \partial_\beta & \partial_\delta 
                      \end{array} \right] \quad;
\end{equation}
\noindent
Their  commutation properties are described by

\begin{equation}\label{dif2}
R_h^{\dagger}Y_2R^{(2) \,-1}Y_1=Y_1R^{(3) \,-1}Y_2R_h \quad,
\end{equation}
\noindent
where $R^{(3)}={\cal P}R^{(2)}{\cal P}$ is given in (\ref{4.8}), 
which is preserved under the $h$-Lorentz coaction; they are  explicitly
given in \cite{BANACH}.

The commutation relations among the entries of $K$ and $Y$ may be expressed 
by an inhomogeneous reflection equation \cite{AKR,FTUV94-21}
\begin{equation}\label{dif4}
Y_2R_hK_1R^{(2)}=R^{(3)}K_1R^{\dagger}_hY_2 +  R^{(3)}{\cal P} \quad,
\end{equation}
\noindent
which extends to the $h$-deformed case the undeformed  relation
$\partial_{\mu}x^{\nu}= \delta^{\nu}_{\mu} + x^{\nu}\partial_{\mu}$.
Eq. (\ref{dif4}) is  consistent with the commutation relations defining the 
algebras  ${\cal M}_h^{(j)}$, ${\cal D}_h^{(j)}$,  and  
is  invariant under $h$-Lorentz transformations. This   
is  seen by multiplying eq. (\ref{dif4})
by $(M_2^{\dagger})^{-1}M_1$  from the left and by $M_1^{\dagger}M_2^{-1}$
from the right and using the commutation relations in (\ref{2.3}).
Eq. (\ref{dif4}) is unique \cite{JPA96} due to the triangularity of $R_h$ 
(in contrast with the situation for $q$-deformation, where covariance allows 
for more than one type of these equations \cite{FTUV94-21,KARPACZ}).
In contrast with the $q$-deformation \cite{OSWZ-CMP,FTUV94-21} (see also
\cite{OZ}), it is a common 
feature of all $h$-deformed   Minkowski spaces that the 
transformation properties for `coordinates' and
`derivatives' are consistent with their simultaneous hermiticity 
\cite{BANACH,PODLES}. 

To determine  the commutation relations for the $h$-de Rham complex we
now introduce the exterior derivative $d$  \cite{OSWZ-CMP,FTUV94-21,AKR}. 
The algebra of  $h$-forms is generated by
the entries of a matrix $dK$. Since $d$ commutes with the 
Lorentz coaction
\begin{equation}\label{dif8}
dK'=M\, dK \, M^{\dagger} \quad.
\end{equation}
\noindent
Applying $d$ to  eq. (\ref{2.6}) we obtain
\begin{equation}\label{dif9}
R_h dK_{1} R^{(2)} K_{2} + R_h K_{1} R^{(2)} dK_{2} = 
dK_{2} R^{(3)} K_{1} R^{\dagger}_h  + K_{2} R^{(3)} dK_{1} R^{\dagger}_h \quad.
\end{equation}
\noindent
Its  only solution is given by
\begin{equation}\label{dif10}
R_h dK_{1} R^{(2)} K_{2} =  K_{2} R^{(3)} dK_{1} R^{\dagger}_h \quad
\end{equation}
\noindent
(which implies  $R_h K_{1} R^{(2)} dK_{2} = 
dK_{2} R^{(3)} K_{1} R^{\dagger}_h$). From eq. (\ref{dif10}) and
$d^2$=0, it follows that 
\begin{equation}\label{dif11}
R_h dK_{1} R^{(2)} dK_{2} = - dK_{2} R^{(3)} dK_{1} R^{\dagger}_h \quad.
\end{equation}
\noindent
Again, it is easy to check that these relations are invariant under 
hermitian conjugation. Notice that eqs. 
(\ref{2.6}), (\ref{dif10})
and (\ref{dif11})  have the same $R$-matrix structure. 
In the $h$-deformed case, the equation giving the commutation 
relations among the  generators of two differential algebras is determined 
only by the transformation (covariant or contravariant) law of these 
generators. Thus, as a consequence  of the triangularity  of $SL_h(2)$,
there are \cite{JPA96,BANACH} only three types of reflection
equations, those in (\ref{2.6}), (\ref{dif2}) and
(\ref{dif4}).
In contrast, this number  is larger for the $SL_q(2)$ based
$q$-deformation.  

For the braiding properties of these algebras we refer to \cite{BANACH}.

\setcounter{equation}{0}

\section{Differential operators and the dilatation operator as a measure 
of the  deformation}

\indent

It is possible to construct from the previous differential algebras some
invariant operators by using the $L_h^{(j)}$-invariant scalar product 
of {\it contravariant} (transforming  as the matrix $K$, eq.  (\ref{2.5}))
and {\it covariant} ($Y \mapsto Y'=(M^{\dagger})^{-1}YM^{-1}$)
matrices (Minkowski four-vectors)  which may be defined as the quantum trace 
\cite{FRT1,ZUM} of a matrix product 
\cite{FTUV94-21}.  The 
$h$-deformed trace \cite{JPA96} of a matrix $B$ is given by 
\begin{equation}\label{gtr}
tr_{h}(B):=tr(D_hB) \quad , \quad  
D_h:=tr_{(2)}({\cal P}(( R_h^{t_1})^{-1})^{t_1})=
\left( \begin{array}{cc} 
        1 & -2h \\
          0   &  1
\end{array} \right) \;, 
\end{equation}
\noindent
where $tr_{(2)}$ means trace in the second space. $tr_h$
is invariant under the quantum group coaction
$B \mapsto MBM^{-1}$ since the expression of $D_h$ above guarantees 
that $D_h^t=M^tD_h^t(M^{-1})^t$ is fulfilled. This is not the only 
possible  definition; for an object $C$ transforming as  
$C \mapsto (M^{\dagger})^{-1}CM^{\dagger}$ we might define another invariant 
trace $\tilde{tr}_h$ by  
\begin{equation}
\tilde{tr}_hC= tr \tilde{D}_hC \quad \mbox{with} \quad 
\tilde{D}_h = \left( \begin{array}{cc} 
        1 & 0 \\
          -2h   &  1
\end{array} \right)  = D_h^{\dagger} \quad.
\end{equation}

The invariant $h$-determinant is related with the $h$-trace \cite{JPA96}. 
Consider
\begin{equation}\label{Reps}
K^{\epsilon}_{ij}:=\hat{R}^{\epsilon}_{h \, ij,kl}K_{kl} \quad, 
\quad  \hat{R}^{\epsilon}_h  :=
(1 \otimes (\epsilon_h^{-1})^t) 
\hat{R}^{(3)} (1 \otimes (\epsilon_h^{-1})^{\dagger}) \quad,
\quad \hat{R}^{\epsilon}_{12}=\hat{R}^{\epsilon}_{21} \;
\end{equation}
\noindent 
(the explicit expressions of $\hat{R}^{\epsilon\,-1}_h$ for (\ref{4.8}) are 
given in eqs. (\ref{hgam1}), (\ref{hgam3})).
Then, if $K$ is contravariant [eq. (\ref{2.5})],  $K^{\epsilon}$ is covariant
{\it i.e.},
$K^{\epsilon} \mapsto (M^{\dagger})^{-1} K^{\epsilon} M^{-1}$.
This may be checked  by using the property of $\hat{R}^{\epsilon}_h$, 
\begin{equation}\label{ap24}
\hat{R}^{\epsilon}_h (M \otimes (M^{\dagger})^t)
=( (M^{\dagger})^{-1} \otimes (M^{-1})^t ) \hat{R}^{\epsilon}_h 
\; \;  \mbox{or} \;\;  \hat{R}^{\epsilon}_h M_1M^*_2=
(M^{\dagger}_1)^{-1} (M^{-1}_2)^t \hat{R}^{\epsilon}_h \;,
\end{equation}
\noindent
which follows from the preservation of the $h$-symplectic metric $\epsilon_h$.
Now, using the expression of  $\epsilon_h$ in (\ref{1.7.4}),
$(P_{h-})_{ij,kl}=- \frac{1}{2}\epsilon_{h\;ij}\epsilon^{-1}_{h\;kl}$ 
and  $D_h= -\epsilon_h(\epsilon_h^{-1})^t$, it
follows that the {\it $h$-deformed  Minkowski length} $l_h$ and 
{\it $h$-metric} $g_h$ are  given by
\begin{equation}\label{hmetr}
l_h := det_hK = \frac{1}{2+h^2}tr_hKK^{\epsilon} \equiv
g_{h\, ij,kl} K_{ij}K_{kl} \quad , \quad
g_{h\, ij,kl}= \frac{1}{2+h^2} D_{h\, si}\hat{R}^{\epsilon}_{h \, js,kl}\;. 
\end{equation}
\noindent
The $h$-metric is preserved under   $h$-Lorentz transformations 
$\Lambda=M \otimes M^*$.

Similarly, the {\it $h$-deformed d'Alembertian} may be  introduced by 
\begin{equation}\label{dif3}
\Box_h := det_hY=  \frac{1}{2+h^2} tr_h(Y^{\epsilon}Y) \quad, \quad 
Y^{\epsilon}= (\hat{R}_h^{\epsilon})^{-1} Y \quad, \quad
Y^{\epsilon} \longmapsto MY^{\epsilon}M^{\dagger} \quad.
\end{equation}
\noindent
As $l_h$,  $\Box_h$  is $h$-Lorentz invariant, real and central 
in the algebra ${\cal D}_h^{(j)}$ of  derivatives. This definition of
$\Box_h$ leads to the Minkowski metric $g^Y$
for the derivatives
\begin{equation}\label{hmetrY1}
\Box_h = g^Y_{h\, ik,mn} Y_{mn}Y_{ik} \quad , \quad
g^Y_{h\, ik,mn}= \frac{1}{2+h^2} D_{h\, kj}
\hat{R}^{\epsilon\,-1}_{h \, ji,mn}\;. 
\end{equation}
\noindent
Using the property $tr_hY^{\epsilon}Y= \tilde{tr}_hYY^{\epsilon}$,
one gets another expression for $g^Y$ which will be used below
\begin{equation}\label{hmetrY2}
g^Y_{h\, ik,mn}= \frac{1}{2+h^2} D_{h\, mj}
\hat{R}^{\epsilon\,-1}_{h \, nj,ik}\;. 
\end{equation}
\noindent
It is easy to check from eqs. (\ref{hmetrY1}) and 
(\ref{hmetrY2}) that  $g^Y$ satisfies $(g^Y)^t={\cal P}g^Y{\cal P}$.
Using eqs. (\ref{hgam1}), (\ref{hgam3}), the $g^{Y\,(j)}$ explicitly read 
\begin{equation}\label{gyexpl}
\begin{array}{c}
\mbox{{\bf 1.}} \;  g^{Y\,(1)}= \frac{1}{2+h^2} \left[
\begin{array}{cccc}
r-h^2 & -h & h & 1  \\
h  & 0 & -1  & 0   \\
-h  & -1 & 0 & 0 \\
1  & 0 & 0 & 0
\end{array} \right] \; ;   \quad
\mbox{{\bf 2.}} \;   
g^{Y\,(2)}=  \frac{1}{2+h^2}  \left[
\begin{array}{cccc}
-5h^2 & 0 & 2h & 1  \\
2h  & 0 & -1  & 0   \\
0  & -1 & 0 & 0 \\
1  & 0 & 0 & 0
\end{array} \right] \;. 
\end{array}
\end{equation}

\vspace{1\baselineskip}

Other invariant differential operators given by the $h$-trace are 
the {\it exterior derivative} $d$, 
\begin{equation}\label{dif12}
d= tr_h(dK\,Y)= d\alpha \partial_\alpha  + d\beta \partial_\beta  +
d\gamma \partial_\gamma  + d\delta \partial_\delta -2h(
d\gamma \partial_\alpha  + d\delta \partial_\beta ) \quad,  
\end{equation}
\noindent
and the {\it dilatation operator} $s$,
\begin{equation}\label{dilat}
s= tr_h(K\,Y)= \alpha \partial_\alpha  + \beta \partial_\beta  +
\gamma \partial_\gamma  + \delta \partial_\delta -2h(
\gamma \partial_\alpha  + \delta \partial_\beta ) \quad.  
\end{equation}

The dilatation operator has a special significance
since it may be considered as a measure of the `strength' of the deformation. 
In the undeformed 
case $s=x^{\mu}\partial_{\mu}$ and accordingly $sx^{\nu}=x^{\nu}(1+s)$.
For the $SL_q(2)$-based  deformations of the Minkowski space, the 
commutation relations of $s_q$ with the generators of ${\cal M}_q$ algebra
acquire extra $q$-dependent terms \cite{FTUV94-21} which do not appear in 
the above undeformed case. For the case of the twisted Minkowski space of 
\cite{CHA-DEM}, the commutation relations of $s$ and $K$ turn out to be 
\cite{FTUV94-21} as in the undeformed case. We shall now prove that this 
behaviour is in fact a general property of the deformed Minkowski algebras
which are defined through a triangular $R$-matrix (this is our case here since
$SL_h(2)$ is triangular \cite{DEMI,ZAKRZ}). 

\vspace{1\baselineskip}

\noindent
{\bf Proposition 3.1}

Let ${\cal M}$ be a deformation of the Minkowski space based on a 
triangular $R$-matrix, and $s$ the dilatation operator defined by 
(\ref{dilat}). Then,
\begin{equation}\label{dilat2}
sK = K(s+1) \quad, \quad  e^{\alpha s}K e^{- \alpha s }=e^{\alpha} K 
\; , \; \alpha \in \bbbr \;.
\end{equation}

\noindent
{\it Proof} : 

\noindent
The invariance of the $h$-trace implies  
\begin{equation}\label{dilat3}
tr_hB=tr_hMBM^{-1}=tr_{h(1)}R_hBR_h^{-1}=tr_{h(1)}R^{(3)}BR^{(3)\,-1} \;
\end{equation}
\noindent
since $R_h$ and $R^{(3)}$ are representations of $GL_h(2)$. Multiplying 
eq. (\ref{dif4}) by $R_{21}K_2$ from the left and by $R^{(2)\,-1}$ from the 
right we get
\begin{equation}\label{dilat4}
R_{21}K_2Y_2R_{12}K_1=R_{21}K_2R^{(3)}K_1R^{\dagger}_{12}Y_2R^{(2)\,-1} 
+  {\cal P}R_{12}K_1 \quad.
\end{equation}
\noindent
We now use eq. (\ref{2.6}) in the first term of the $r.h.s.$  and take the 
$h$-trace in the second space. Using that for a triangular $R$-matrix 
$R_{12}=R_{21}^{-1}$ the $tr_{h(2)}$ satisfies $tr_{h(2)}R_{21}BR_{12}
=tr_{h(2)}B$, we obtain
\begin{equation}\label{dilat5}
s \, K_1 = K_1 \, tr_{h(2)} R^{(2)}K_2Y_2R^{(2)\,-1} + 
K_1\, tr_{h(2)}{\cal P}R_{12} \quad.
\end{equation}
\noindent
Since  $tr_{h(2)}{\cal P}R_{12}\equiv tr_{(2)}(I_2 \otimes D_h)\hat{R}_h=I_2$
(for $R_q$ this is $q^2I_2$), we obtain $sK=K(1+s)$, {\it q.e.d.}.

\setcounter{equation}{0}

\section{$h$-Minkowski algebras and representations}
 
\indent

Using the $R^{(3)}$ matrices given in (\ref{4.8}) in  eq. (\ref{2.6}),
the Minkowski algebras associated with $SL_h(2)$ explicitly read 

\noindent
{\bf 1.} $\quad {\cal M}_h^{(1)}$:
\begin{equation}\label{5.7}
\begin{array}{ll}
{} [ \alpha, \beta ]= -h \beta^2 -r \beta \delta + h \delta \alpha 
-h \beta \gamma + h^2 \delta \gamma \;,&
\quad [\alpha, \delta]= h( \delta \gamma - \beta \delta) \;,\\
{} [\alpha, \gamma]= h \gamma^2 + r \delta \gamma - h \alpha \delta
+h \beta \gamma - h^2 \beta \delta \;, & \quad 
[\beta, \delta]=  h \delta^2 \;,\\
{} [ \beta, \gamma]= h \delta ( \gamma + \beta) + r \delta^2 \;,& \quad
[ \gamma, \delta]= - h \delta^2 \;;
\end{array}
\end{equation}

\vspace{1\baselineskip}

\noindent
{\bf 2.} $\quad {\cal M}_h^{(2)}$:
\begin{equation}\label{5.10}
\begin{array}{ll}
{} [ \alpha, \beta ]= 2h \alpha \delta + h^2 \beta \delta  \quad , &
\quad [\alpha, \delta]= 2h (\delta \gamma -  \beta \delta ) \quad,\\
{} [\alpha, \gamma]= -h^2 \delta \gamma - 2h \delta \alpha \quad, & \quad 
[\beta, \delta]=  2 h \delta^2 \quad , \\
{} [ \beta, \gamma]= 3h^2 \delta^2  \quad ,& \quad
[ \gamma, \delta]= - 2h \delta^2 \quad .
\end{array}
\end{equation}

\noindent
In each case there exists a subalgebra generated by $\beta,\gamma,\delta$
which will be denoted by ${\cal A}^{(j)}$ ($j$=1,2).

The $h$-deformed Minkowski length determines quadratic central elements 
for both  ${\cal M}_h^{(1)}$ and ${\cal M}_h^{(2)}$,
\begin{equation}\label{5.7.1}
l_h^{(1)}= \frac{2}{h^2+2}(\alpha \delta - \beta \gamma + h \beta \delta)
\quad,
\end{equation}
\begin{equation}\label{5.11}
l_h^{(2)}=\frac{2}{h^2+2}(\alpha \delta - \beta \gamma + 2h \beta \delta) 
\quad.
\end{equation}

For ${\cal M}_h^{(2)}$, it is possible to find a linear central real element 
given by
\begin{equation}\label{lin}
\zeta := \beta + \gamma - \frac{3}{2}h \delta \quad.
\end{equation}
\noindent
We may extract from ${\cal A}^{(2)}$ the following set of hermitian 
generators 
\begin{equation}\label{gen1}
x:=h \delta \quad , \quad y:= -i(\beta-\gamma) \quad,\quad 
\zeta := \beta + \gamma - \frac{3}{2}h \delta \quad,
\end{equation}
\noindent
($\alpha^* = \alpha$,  $\delta^* =\delta$, $\beta^*=\gamma$, eq. (\ref{2.5}))
and hence
\begin{equation}\label{gen2}
\delta=x/h \quad,\quad \beta= \frac{1}{2} (\frac{3}{2}x+iy+ \zeta) \quad, \quad
\gamma = \frac{1}{2} (\frac{3}{2}x-iy+ \zeta) \;.
\end{equation}
\noindent
This makes easier to discuss the representations of ${\cal M}_h^{(2)}$. Since 
$\zeta$ is central, we can represent it by a real
number $\zeta \in \bbbr$. The commutation relations of $x$, $x^{-1}$ and $y$ 
are
\begin{equation}\label{xy1} 
[x,y]=4ix^2 \quad,\quad [y,x^{-1}]=4i \quad.
\end{equation}
\noindent
The second commutator is of the  canonical Heisenberg  type, so the unique
(up to unitary equivalence) irreducible representation is $L_2(\bbbr)$
(Heisenberg algebra) with $x$ as the multiplication operator and 
\begin{equation}\label{xy2} 
y=-4ix \frac{d}{dx} x = -4i(x^2 \frac{d}{dx} + x )\quad.
\end{equation}
\noindent
We may now extend this representation of ${\cal A}^{(2)}$
to the whole algebra ${\cal M}_h^{(2)}$. The generator $\alpha$  in
this irreducible representation is defined by fixing  the quadratic central 
element $l_h^{(2)} \in \bbbr$ (if we consider ${\cal M}_h$ as a
$h$-deformed momenta algebra this condition corresponds to fixing the 
mass-shell). Using
(\ref{gen1}) in the expression of $l_h^{(2)}$ [eq. (\ref{5.11})] one finds
\begin{equation}\label{xy3}
\begin{array}{c}
x \alpha = h (\frac{2+h^2}{2})l_h^{(2)} + h (\beta -2x) \gamma
\end{array}
\end{equation}
\noindent
and, using (\ref{gen2}) to express $\beta$ and $\gamma$ in terms of  $x$,
$y$ and $\zeta$
\begin{equation}\label{xy4}
\alpha = \frac{h}{4}( \frac{1}{x}(y^2 + 2i \{ y,x \}_+)- \frac{23}{4}x- \zeta
+ \frac{1}{x}( \zeta^2 + (4+h^2) l_h^{(2)}))  \quad.
\end{equation}
\noindent
There is also a one-dimensional representation with
\begin{equation}\label{xy5}
\alpha \in \bbbr \quad , \quad \delta =0 \quad , \quad  \beta = \gamma^*
\, \in \bbbc \;.
\end{equation}
\noindent
These representations are very different from those of $q$-deformed Minkowski
spaces  \cite{FTUV94-21}, but also from other twisted Minkowski
spaces as that of \cite{CHA-DEM}  (called ${\cal M}_p$ in \cite{FTUV94-21}).

\setcounter{equation}{0}

\section{Properties of the $h$-Lorentz algebras: triangularity and twisted 
character}

\indent

  In the earliest papers \cite{WA-ZPC48,OSWZ-CMP} on the 
deformed spacetime algebra 
${\cal M}_q$ (associated with $SL_q(2)$) the commutation relation of 
the quantum coordinates (generators of ${\cal M}_q$) were written 
in an exchange algebra form
\begin{equation}\label{ex1}
{\cal R}_{12}X_1X_2= \mu X_2X_1
\end{equation}
\noindent
in terms of a 16$\times$16 $R$-matrix ${\cal R}$
\cite{OSWZ-CMP} (or the corresponding projector \cite{WA-ZPC48})
and a four component column vector $X$. The convenience of the reflection 
equation formalism (see {\it e.g.} \cite{KS}) to achieve an unified 
description of the different deformations of the Lorentz group and the 
Minkowski space for dimensionless deformation parameters was exhibited
in a series of papers \cite{AKR,FTUV94-21,JPA96,PLB95}. The commutation 
relations of the deformed coordinates were written  in the generic form
\begin{equation}\label{ex2}
R^{(1)} K_{1} R^{(2)} K_{2} = K_{2} R^{(3)} K_{1} R^{(4)} \quad,
\end{equation}
\noindent
where  the 4$\times$4 matrices $R^{(i)}$ ($i$=1,2,3,4) satisfy a number 
of consistency conditions  (as (\ref{2.4})) \cite{FM,FTUV94-21,JPA96},
and $K$ is a 2$\times$2 matrix  representing the deformed coordinates 
(eq. (\ref{2.6}) is a particular case of this procedure).

In order to have an explicit comparison of this approach with the above 
mentioned papers \cite{WA-ZPC48,OSWZ-CMP} as well as to study the 
triangularity of $L_h^{(j)}$, it is 
convenient to transform the reflection equation (\ref{ex2}) into the
exchange algebra (\ref{ex1})  and to relate the properties of the $R$-matrix
${\cal R}$ with those of the $R^{(i)}$.
Provided the existence of the corresponding inverse matrices, eq. (\ref{ex2})
\begin{equation}\label{ex3}
R^{(1)}_{ij,ab} K_{ac} R^{(2)}_{cb,kd} K_{dl} =  
K_{ja'} R^{(3)}_{ia',b'c'} K_{b'd'} R^{(4)}_{d'c',kl} \quad
\end{equation}
\noindent
(summation over repeated indices  understood) may be written in the exchange
algebra form
\begin{equation}\label{ex4}
K_{ac} K_{dl} =  {\cal R}_{(dl)(ac),(ja')(b'd')} K_{ja'} K_{b'd'} \quad,
\end{equation}
\noindent
where
\begin{equation}\label{ex5}
{\cal R}_{(dl)(ac),(ja')(b'd')}=
(R^{(1)})^{-1}_{ab,ij} ( R^{(2)\,t_2})^{-1}_{kb,cd} 
 R^{(3)}_{ia',b'c'}  R^{(4)}_{d'c',kl} \quad
\end{equation}
\noindent
(see \cite{UMeyer} for an analogous discussion in the particular case of 
the $q$-Minkowski space of \cite{WA-ZPC48,OSWZ-CMP}). Eq. (\ref{ex4}) 
corresponds 
to the form  (\ref{ex1}) once the matrix elements
$K_{ab}$ are understood as components of a column vector labelled by a pair
of indices [$(ab)=(11),(12),(21),(22)$].
Due to the various consistency conditions \cite{FTUV94-21,FM} that the 
$R^{(i)}$ matrices have to satisfy, the 
${\cal R}$ satisfies the Yang-Baxter equation.

\vspace{1\baselineskip}

\noindent
{\it Remark}. In the case where all the $R^{(j)}$ ($j$=1,2,3,4) matrices in
(\ref{ex2}) are expressed in terms of the standard $R_q$ matrix \cite{FRT1}
for $GL_q(2)$, and which leads to the Minkowski algebra ${\cal M}_q$
\cite{WA-ZPC48,OSWZ-CMP}, considerations on the projectors \cite{OSWZ-CMP}
lead to the introduction of two different 16$\times$16 $R$-matrices 
associated with the $q$-Lorentz group. In this case,
the matrix ${\cal R}$ of (\ref{ex5}) coincides with one of them
(after a change of basis it corresponds to $R_{II}$ in \cite{OSWZ-CMP}).
To find the other $R$-matrix  (which is the one used for the definition of the 
$q$-Lorentz group in \cite{WA-ZPC48}) the starting point is not (\ref{ex2}),
but the reflection equation which provides the commutation relations 
between the differentials $dK$ (\cite{AKR}, eq. (37)). Nevertheless, 
the procedure
leading to (\ref{ex5}) is the same  (for a comparison of the eqs. in
\cite{OSWZ-CMP} with this formalism see \cite{KARPACZ}).

\vspace{1\baselineskip}

The triangularity property of the $R$-matrix, $R_{12}R_{21}=I$ \cite{DRINF},
is important in the classification scheme of deformed Poincar\'e algebras
\cite{PODWOR} and has also been used in a recent proposal to solve  
the Klein-Gordon and 
Dirac equations in deformed Minkowski spaces \cite{PODLES}. The
description of ${\cal M}_h^{(j)}$ of Sec.2  related to the `triangular'
quantum group $SL_h(2)$ is, on the other hand, based on the matrices 
$R^{(j)}$, of which $R^{(1)}=R_h$ and $R^{(4)}$ are triangular. It is
then natural to ask whether the triangularity of the 16$\times$16
matrix ${\cal R}$ is implied by the triangularity of the 4$\times$4 
$R$-matrices. Due to the general consistency relation $R^{(3)}
={\cal P}R^{(2)}{\cal P}$  the following proposition holds:

\vspace{1\baselineskip}

\noindent
{\bf Proposition 5.1} 

 Let $R^{(1)}$ and $R^{(4)}$ be triangular $R$-matrices  ($R_{12}R_{21}=I$)
and let $R^{(3)}$ and $R^{(2)}$ related by 
$R^{(3)} = {\cal P}R^{(2)}{\cal P}$ and satisfy the consistency conditions 
(\ref{2.4}) (eq.  (\ref{2.4}) is not necessary below, but it is implicit in 
the introduction of ${\cal R}$). Then, the $R$-matrix ${\cal R}$ 
in (\ref{ex5}) is also triangular {\it i.e.}, 
\begin{equation}\label{tri1}
{\cal R}_{(dl)(ac),(ja')(b'd')} {\cal R}_{(b'd')(ja'),(mn)(rs)}  =
 \delta_{dr} \delta_{ls} \delta_{am} \delta_{cn}=I_{(dl)(ac),(rs)(mn)}
\;
\end{equation}
\noindent
since $({\cal R}_{21})_{(ac)(dl),(b'd')(ja')}=
({\cal R}_{12})_{(dl)(ac),(ja')(b'd')}$.

\noindent  
{\it Proof}: 

\noindent
Using (\ref{ex5}), the ${\cal R}$ matrices in the $l.h.s.$
of (\ref{tri1}) are expressed in terms of $R^{(i)}$ ($i$=1,2,3,4)
by
\begin{equation}\label{tri2}
(R^{(1)})^{-1}_{ab,ij} ( R^{(2)\,t_2})^{-1}_{kb,cd} 
 R^{(3)}_{ia',b'c'}  R^{(4)}_{d'c',kl} 
(R^{(1)})^{-1}_{jp,tm} ( R^{(2)\,t_2})^{-1}_{zp,a'b'} 
R^{(3)}_{tn,rf}  R^{(4)}_{sf,zd'} \;.
\end{equation}
\noindent
Thus, the proof of (\ref{tri1}) reduces to some index contractions between
pairs of $R$-matrices. We start by 
\begin{equation}\label{tri3}
 R^{(3)}_{ia',b'c'}( R^{(2)\,t_2})^{-1}_{zp,a'b'} =
 R^{(2)\,t_2}_{a'b',c'i}( R^{(2)\,t_2})^{-1}_{zp,a'b'}
=\delta_{ip}\delta_{zc'} \;.
\end{equation}
\noindent
Then we have independent contractions in which we must use the triangularity
of $R^{(1)}$ and $R^{(4)}$ 
\begin{equation}\label{tri4}
(R^{(1)})^{-1}_{ab,ij}(R^{(1)})^{-1}_{ji,tm}=
(R^{(1)})^{-1}_{ab,ij} R^{(1)}_{ij,mt}= \delta_{am} \delta_{bt} \quad,
\end{equation}
\begin{equation}\label{tri5}
R^{(4)}_{d'c',kl}  R^{(4)}_{sf,c'd'}=R^{(4)}_{d'c',kl}
(R^{(4)})^{-1}_{fs,d'c'}= \delta_{kf} \delta_{ls} \quad.
\end{equation}
\noindent
Finally, the last contraction is similar to the first one
\begin{equation}\label{tri6}
( R^{(2)\,t_2})^{-1}_{kb,cd}  R^{(3)}_{bn,rk} =
( R^{(2)\,t_2})^{-1}_{kb,cd}  R^{(2)\,t_2}_{nr,kb}= \delta_{cn} 
\delta_{dr} \;.
\end{equation}
\noindent
{\it q.e.d.}

\vspace{1\baselineskip}

\noindent
{\it Remark}. We stress that the proof depends only on the triangularity   of
$R^{(1)}$ and $R^{(4)}$ (which, in a more general setting than the present 
one in  which we are concerned with deformed spacetimes, might even
have different dimensions) and not $R^{(2)}$ 
($R^{(3)}={\cal P}R^{(2)}{\cal P}$).

\vspace{1\baselineskip}

Using the expression of $R_h=R^{(1)}=R^{(4)\,t}$  (eq. (\ref{1.5}))
and $R^{(3)\,(j)}$ ($j$=1,2)  (eq. (\ref{4.8})) in (\ref{ex5})
the explicit form of ${\cal R}$ is found. The two matrices ${\cal R}^{(j)}$ 
for the $h$-deformations $L_h^{(j)}$  are given in 
the  appendix, where  their triangularity may be checked by 
direct computation.

It is well known \cite{DRINF,RESH} that a 
special class of triangular $R$-matrices is obtained by the twisting
procedure. A twisted $R$-matrix is obtained from a 
non-singular matrix $F$ (with some extra properties) by
\begin{equation}\label{tw1}
R_{12}:=F{\cal P}F^{-1}{\cal P}= F \tilde{F}^{-1} \quad , \quad 
\tilde{F}:={\cal P}F{\cal P} \quad.
\end{equation}
\noindent
This is the case of the $R$-matrix of $SL_h(2)$ \cite{DEMI}.
The $R$-matrix $R_h$ (\ref{1.5}) may be expressed in terms of a matrix
$F$ which is given in the fundamental representation by \cite{ZAKRZ}
\begin{equation}\label{tw2}
F= \left[ 
\begin{array}{cccc}
1 & h/2 & -h/2 & 0 \\
0 & 1 & 0 & h/2 \\
0 & 0 & 1 & -h/2 \\
0 & 0 & 0 & 1 
\end{array} \right] = I+ \frac{h}{2}(H \otimes E - E \otimes H) \;
\end{equation}
\noindent
where $E= \left[ 
\begin{array}{cc}
0 & 1 \\
0 & 0
\end{array} \right]$ and  $H = \left[ 
\begin{array}{cc}
1 & 0 \\
0 & -1
\end{array} \right]\,$.

Since for the $h$-Minkowski algebras $R^{(4)}=R^{t}_h$ (see (\ref{2.6}) and
(\ref{ex2})) it follows that $R^{(4)}$ is also obtained by twisting from a 
matrix $G$ related to $F$, $R^{(4)}=R^{(1)\,t}= 
\tilde{F}^{-1\,t}F^t = G \tilde{G}^{-1}$ so that
\begin{equation}\label{twg} 
G=({\cal P}F^{t}{\cal P})^{-1}=
\left[ 
\begin{array}{cccc}
1 & 0 & 0 & 0 \\
h/2 & 1 & 0 & 0 \\
-h/2 & 0 & 1 & 0 \\
h^2/2 & h/2 & -h/2 & 1 
\end{array} \right] \quad.
\end{equation}

Let us now prove that, much in the same way as the triangularity property
(Proposition 5.1), the property that $R^{(1)}$ and $R^{(4)}$ are obtained
by twisting is also inherited by the 16$\times$16 matrix ${\cal R}$ of
the deformed Lorentz group. This property, as the previous triangularity one,
is general when the $R^{(1)}$ and $R^{(4)}$ are those appearing in a 
reflection equation  (cf. (\ref{ex2})) and ${\cal R}$ is the matrix 
appearing in the associated  exchange algebra.

\vspace{1\baselineskip}

\noindent
{\bf Proposition 5.2}

Let $R^{(1)}$ and $R^{(4)}$  in (\ref{ex2}) be obtained by a twisting 
procedure, {\it i.e.}, there exist two invertible matrices $F$ and $G$
such that
\begin{equation}\label{tw3}
R^{(1)}= F \tilde{F}^{-1} \quad , \quad R^{(4)}= G \tilde{G}^{-1} 
\quad.
\end{equation}
\noindent
Then the exchange algebra 16$\times$16
matrix ${\cal R}$ given by (\ref{ex5}) has the same structure 
\begin{equation}\label{tw4}
{\cal R}= {\cal F} \tilde{{\cal F}}^{-1} \quad 
\end{equation}
\noindent
where
\begin{equation}\label{tw5}
{\cal F}_{(dl)(ac),(rm)(sn)}= F_{ba,sr} (R^{(2)\,t_2})^{-1}_{kb,cd}
\tilde{G}^{-1}_{mn,kl} \quad.  
\end{equation}
\noindent
{\it Proof} :  

\noindent
Using (\ref{tw3})  in (\ref{ex5}) we obtain an expression
of ${\cal R}$ as a product of two matrices 
\begin{equation}\label{tw6}
\begin{array}{ll}
{\cal R}_{(dl)(ac),(ja')(b'd')} & =
\tilde{F}_{ab,rs}F^{-1}_{rs,ij} ( R^{(2)\,t_2})^{-1}_{kb,cd} 
 R^{(3)}_{ia',b'c'}  G_{d'c',mn} \tilde{G}^{-1}_{mn,kl} \\
\, &  =[\tilde{F}_{ab,rs} ( R^{(2)\,t_2})^{-1}_{kb,cd} \tilde{G}^{-1}_{mn,kl}]
[F^{-1}_{rs,ij} R^{(2)\,t_2}_{a'b',c'i} G_{d'c',mn}] \\
\, & := {\cal F}_{(dl)(ac),(rm)(sn)}{\cal H}_{(rm)(sn),(ja')(b'd')} \quad.
\end{array} 
\end{equation}
\noindent
To prove the twisted character of ${\cal R}$ we must now check that 
${\cal H}=\tilde{{\cal F}}^{-1}$. Multiplying $\tilde{{\cal F}}$ (obtained
from (\ref{tw5}), $\tilde{{\cal F}}_{(ac)(dl),(sn)(rm)} 
={\cal F}_{(dl)(ac),(rm)(sn)}$) and ${\cal H}$, one gets
\begin{equation}\label{tw7}
\begin{array}{l}
\tilde{{\cal F}}_{(dl)(ac),(rm)(sn)}{\cal H}_{(rm)(sn),(ja')(b'd')}  
 =  \tilde{F}_{bd,sr} ( R^{(2)\,t_2})^{-1}_{kb,la} \tilde{G}^{-1}_{nm,kc}
F^{-1}_{rs,ij} R^{(2)\,t_2}_{a'b',c'i} G_{d'c',mn} \\
\qquad = \delta_{bi}\delta_{dj}\delta_{c'k}\delta_{d'c} 
( R^{(2)\,t_2})^{-1}_{kb,la} R^{(2)\,t_2}_{a'b',c'i} 
 =   \delta_{dj}\delta_{d'c} 
( R^{(2)\,t_2})^{-1}_{c'i,la} R^{(2)\,t_2}_{a'b',c'i} \\
\qquad =   \delta_{dj}\delta_{d'c}\delta_{a'l}\delta_{b'a} =
I_{(dl)(ac),(ja')(b'd')} \quad,
\end{array}
\end{equation}
\noindent
{\it q.e.d.}
 
\vspace{1\baselineskip}

\noindent
{\bf Corollary 5.1} \,
The $h$-deformations of the Lorentz group associated with $SL_h(2)$ 
\cite{WOZA,JPA96} are a twisting of the standard Lorentz group.

\setcounter{equation}{0}

\section{Dirac $\gamma$-matrices for $h$-Minkowski spaces}

\indent

Any discussion of deformed relativistic equations requires the expression 
of the deformed d'Alembertian for the Klein-Gordon equation and the deformed 
Dirac matrices for the Dirac one. In the present case, $\Box_h$ is given 
by eq. (\ref{dif3}). We shall devote this section to find the 
explicit form of the $h$-Dirac matrices and to prove that they satisfy
the appropriate $h$-anticommutation properties.

   The generators of a  deformed Minkowski algebra ${\cal M}_h^{(i)}$
 and the generators of the corresponding derivatives 
algebra ${\cal D}_h^{(i)}$ can be arranged in the elements of  2$\times$2 
matrices $K$ and $Y$ of (\ref{2.6}) and (\ref{dif2}).
Using the set of 2$\times$2 matrices  $e_{ij}$ ($i,j$=1,2) defined by
$(e_{ij})_{\alpha \beta}=\delta_{i \alpha} \delta_{j \beta}$, as a basis
of the 2$\times $2 matrices, we may write
\begin{equation}\label{gam1}
Y=e_{ij} \partial_{ij} \quad.
\end{equation}
\noindent
Any linear transformation of the generators $\partial_{ij} \rightarrow 
\partial'_{mn}= A_{mn,ij} \partial_{ij}$ results in 
the inverse transformation  for the matrices $e_{ij}$.

The quantum group coaction for a covariant vector $Y \rightarrow 
(M^{\dagger})^{-1}YM^{-1}$ (with $det_hY \neq 0$) obviously gives for 
$Y^{-1}$ the transformation law of a contravariant vector $Y^{-1} \rightarrow
MY^{-1}M^{\dagger}$. A contravariant vector $Y^{\epsilon}$ may be constructed
as a linear combination of the entries of $Y$ by
\begin{equation}\label{gam2}
Y^{\epsilon}_{ij}=(\hat{R}^{\epsilon}_h)^{-1}_{ij,kl}Y_{kl} \quad,
\end{equation}
\noindent
(cf. (\ref{Reps})).  It may be seen that 
the contravariant vectors $Y^{-1}$ and $Y^{\epsilon}$ are related by an
invariant and central factor,
\begin{equation}\label{yye} 
Y^{\epsilon}= \frac{2+h^2}{2} (det_hY)Y^{-1}= \frac{2+h^2}{2} \Box_h Y^{-1} 
\quad.
\end{equation}

A natural definition of the deformed Dirac operator (see \cite{FTUV94-21})
is the following 
\begin{equation}\label{gam3}
\Dslash = \left[  \begin{array}{cc} 0 & Y \\ Y^{\epsilon} & 0 
\end{array} 
\right]  \quad.
\end{equation}
\noindent
This definition is in fact general and is also valid for all 
the spacetime $Q$-deformations  ($Q$=$q,h$) in \cite{JPA96};
the different cases differ from one another in the explicit expression of 
$Y^{\epsilon}$ and on the commutation relations among the generators of the 
${\cal D}_Q$  algebras obtained for the different deformed $Y$'s.

Using eqs. (\ref{gam1}) and (\ref{gam2}) we find that the 2$\times$2 matrix
$Y^{\epsilon}$  is given by
\begin{equation}\label{gam4}
Y^{\epsilon}= [e_{kl}(\hat{R}^{\epsilon}_{h})^{-1}_{kl,ij}] \partial_{ij} \quad.
\end{equation}
\noindent
Thus, the Dirac operator can be written as a contraction of the initial 
generators $\partial_{ij}$ with 4$\times$4 matrices (deformed $\gamma$-matrices)
\begin{equation}\label{gam5}
\Dslash  = \gamma_{ij} \partial_{ij} \quad, \quad
\gamma_{ij} := \left( \begin{array}{cc} 0 & 1 \\ 0 & 0 \end{array} \right)
\otimes e_{ij} + \left( \begin{array}{cc} 0 & 0 \\ 1 & 0 \end{array} \right)
 \otimes e_{kl}(\hat{R}^{\epsilon}_{Q})^{-1}_{kl,ij} \;.
\end{equation}
\noindent
Clearly, the more familiar notation is obtained by replacing the 
pairs $(ij),(mn)$... by $\mu , \nu$...

We have to  prove  that the $h$-deformed $\gamma$-matrices satisfy 
the suitable
deformed Clifford algebra relation; as in the undeformed case, it must be 
consistent 
with the relation $\Dslash^2 \propto \Box_h$ and with the commutation 
relation for the derivatives. Hence, the deformed relation for 
$\gamma$-matrices  will be written in terms of the metric given for $Y$, 
eq. (\ref{hmetrY1}), and 
in terms of the 16$\times$16 ${\cal R}$-matrix
of the exchange algebra for $Y$, ${\cal R}^Y$. This is obtained from 
the reflection equation  for $Y$ (see  (\ref{dif2})) 
\begin{equation}\label{gam6}
R^{(4)}Y_2R^{(2) \,-1}Y_1=Y_1R^{(3) \,-1}Y_2R^{(1)} \quad,
\end{equation}
\noindent
($R^{(3)}={\cal P}R^{(2)}{\cal P}$) analogously to (\ref{ex5}), by 
identifying the intermediate expression 
\begin{equation}
Y_{ia'}Y_{c'd'} R^{(1)}_{b'd',mn}((R^{(2) \,-1})^{t_2})^{-1}_{dc,an} R^{(3) \,-1}_{a'j,b'c'}
R^{(4)\,-1}_{ab,ij} =  Y_{bc}Y_{dm}
\end{equation}
\noindent
with  the exchange algebra relation $Y_1Y_2 {\cal R}^Y = Y_2Y_1$. Then,
\begin{equation}\label{gam7}
{\cal R}^Y_{(ia')(c'd'),(dm)(bc)} =
R^{(1)}_{b'd',mn}((R^{(2) \,-1})^{t_2})^{-1}_{dc,an} R^{(3) \,-1}_{a'j,b'c'}
R^{(4)\,-1}_{ab,ij} \quad.
\end{equation}

\vspace{1\baselineskip}

\noindent
{\bf Proposition 6.1}

The $h$-deformed Dirac matrices  (\ref{gam5}) satisfy the deformed Clifford 
algebra relations
\begin{equation}\label{clif1}
\gamma_{ij}\gamma_{mn} + {\cal R}^Y_{(ij)(mn),(rs)(kl)}\gamma_{kl}
\gamma_{rs} = (2+h^2) g^Y_{h \, mn,ij} I_4 \quad,
\end{equation}
\noindent
where ${\cal R}^Y$ and $g^Y_h$ are given by (\ref{gam7}) and (\ref{hmetrY2})
respectively.

\noindent
{\it Proof} :  

\noindent
Using the Weyl-like realization (\ref{gam5}) for the deformed $\gamma$-matrices,
eq. (\ref{clif1}) splits into two equations involving the 2$\times$2 
matrices $e_{ij}$, namely
\begin{equation}\label{clif2}
e_{ij}e_{ab} (\hat{R}^{\epsilon \, -1}_h)_{ab,mn} + 
{\cal R}^Y_{(ij)(mn),(rs)(kl)}e_{kl}e_{cd} (\hat{R}^{\epsilon \, -1}_h)_{cd,rs} 
= (2+h^2) g^Y_{h \, mn,ij} I_2 \quad,
\end{equation}
\begin{equation}\label{clif3}
e_{zt} (\hat{R}^{\epsilon \, -1}_h)_{zt,ij} e_{mn}+ 
{\cal R}^Y_{(ij)(mn),(rs)(kl)}e_{fg} (\hat{R}^{\epsilon \, -1}_h)_{fg,kl}e_{rs} 
= (2+h^2) g^Y_{h \, mn,ij} I_2 \quad,
\end{equation}
\noindent 
Let us prove eq. (\ref{clif2}) (eq. (\ref{clif3}) is proved similarly).
Eq. (\ref{clif2}) 
may be written in components by using the explicit form of the matrices
$(e_{ij})_{\alpha \beta}= \delta_{i \alpha} \delta_{j \beta}$ so that
$(e_{ij}e_{ab})_{\alpha \beta}=
(e_{ij})_{\alpha \gamma}(e_{ab})_{\gamma \beta}=
 \delta_{i \alpha} \delta_{ja} \delta_{b \beta}$. Then, the matrix element 
$\alpha \beta$ of the $l.h.s.$ of the eq.  (\ref{clif2}) is
\begin{equation}\label{clif4}
\delta_{i \alpha} (\hat{R}^{\epsilon \, -1}_h)_{j \beta,mn} + 
{\cal R}^Y_{(ij)(mn),(rs)(\alpha c)} (\hat{R}^{\epsilon \, -1}_h)_{c \beta,rs} 
\quad .
\end{equation}
\noindent
Now, using the expression (\ref{gam7}) for ${\cal R}^Y$ and 
(\ref{Reps}) for $\hat{R}^{\epsilon}_h$, one gets for the $l.h.s.$
\begin{equation}\label{clif5}
\delta_{i \alpha} \epsilon^t_{h \, \beta \beta'} R^{(3)\,-1}_{ j \beta', n'm}
\epsilon^t_{h \, n'n} + R^{(1)}_{b'n,sl}((R^{(2) \,-1})^{t_2})^{-1}_{rc,al} 
R^{(3) \,-1}_{jk,b'm} R^{(4)\,-1}_{a \alpha ,ik} 
\epsilon^t_{h \, \beta \beta'} R^{(3)\,-1}_{ c \beta', s'r}
\epsilon^t_{h \, s's} \quad.
\end{equation}
\noindent
contracting pairs of indices in the second term of (\ref{clif5})
\begin{equation}\label{clif6}
((R^{(2) \,-1})^{t_2})^{-1}_{rc,al} R^{(3)\,-1}_{ c \beta', s'r} =
((R^{(2) \,-1})^{t_2})^{-1}_{rc,al} (R^{(2)\,-1})^{t_2}_{\beta' s',rc} =
\delta_{\beta' a } \delta_{s' l} \quad, 
\end{equation}
\noindent
and using
\begin{equation}\label{clif7}
R^{(1)}_{b'n,sl} \epsilon_{h \, sl} = -\epsilon_{h \, nb'} \quad, 
\end{equation}
\noindent
(which follows from (\ref{1.7.2}), (\ref{1.7.4})) allows us 
to write the $l.h.s.$  of  (\ref{clif2}) as 
\begin{equation}\label{clif8}
\epsilon^t_{h \, \beta \beta'} \left( \delta_{i \alpha} \delta_{\beta' k}
- R^{(4)\,-1}_{\beta' \alpha, ik} \right) R^{(3)\,-1}_{jk,b'm}
\epsilon_{h \, nb'} \quad.
\end{equation}
\noindent
Since $R^{(4)}=R^{(1)\, \dagger}=R_h^t$, the term in the brackets
is $(I- \hat{R}_h^{-1})_{ik, \alpha \beta'}$  which,  using the spectral 
decomposition of $\hat{R}_h$ (eq. (\ref{1.7.2})), is equal to
$2P_{h-\,ik, \alpha \beta'}$. Finally, using eq. (\ref{1.7.4}) to express
$P_{h-}$ in terms of $\epsilon_h$, one gets by direct computation
(see eq. (\ref{hmetrY2}))
\begin{equation}\label{clif9}
\delta_{\alpha \beta} \left( - \epsilon_{h\,iz}\epsilon_{h\,sz}^{-1} \right)
\left( \epsilon_{h\,sk}^t  \hat{R}^{(3)\,-1}_{jk,mb'}
 \epsilon_{h\,b'n}^t \right) = \delta_{\alpha \beta} D_{h\,is}
(\hat{R}^{\epsilon\,-1}_h)_{js,mn} = \delta_{\alpha \beta} (2+h^2)
g^Y_{h\, mn,ij} 
\end{equation}
\noindent
which coincides with the matrix element $\alpha \beta$ of the $r.h.s.$ of 
eq. (\ref{clif2}), {\it q.e.d.}

\vspace{1\baselineskip}

If instead of the basis $ \{ e_{ij} \}$  
another one is selected, the resulting matrices will be the gamma matrices
in another `spacetime' basis. However, in the deformed case 
there is not a clear way
to define a `physical' basis. For the
$q$-Minkowski space of \cite{WA-ZPC48,OSWZ-CMP}, there is a central element 
in the algebra  which may be associated with time and the other generators 
may be grouped by its tensorial properties under the deformed rotation
subgroup, so that one may define $q$-Pauli and $q$-gamma matrices  
(see in this respect \cite{WA-ZPC48,SONG,SCHI,FTUV94-21}). However, 
there is no
deformed rotation subgroup for all Lorentz deformations (as {\it e.g.},
for \cite{CHA-DEM}), and there is no central element that we can clearly select
as time, as in the present case of the $h$-deformation.  Thus, we shall not 
attempt to introduce any special basis and we shall work with the generators
$\alpha$, $\beta$, $\gamma$, $\delta$ in (\ref{2.5}). Thus, we shall label 
the four gamma matrices as $\gamma^{\alpha}$, $\gamma^{\beta}$,
$\gamma^{\gamma}$, $\gamma^{\delta}$ and write $\Dslash =\gamma^I \partial_I$
($I$ = $\alpha$, $\beta$, $\gamma$, $\delta$). 

To find their explicit form, we need the  expressions of 
$\hat{R}^{\epsilon}_h$ and $Y^{\epsilon}$, which are found from (\ref{Reps}) 
and (\ref{gam2}) for ${\cal M}_h^{(1)}$ and ${\cal M}_h^{(2)}$. The 
results are

\vspace{1\baselineskip}

\noindent
{\bf 1. \,}  ${\cal M}_h^{(1)}$  case:

\begin{equation}\label{hgam1}
(\hat{R}_h^{\epsilon})^{-1} = \left[ \begin{array}{cccc}
h^2+r & -h & -h & 1 \\
h & -1 & 0 & 0 \\
h & 0 & -1 & 0 \\
1 & 0 & 0 & 0  \end{array} \right] \; ; \;
Y^{\epsilon} = \left[ \begin{array}{cc}
(h^2+r) \partial_{\alpha}  -h(\partial_{\beta}+ \partial_{\gamma}) 
+ \partial_{\delta} & h \partial_{\alpha}  - \partial_{\gamma}  \\
h \partial_{\alpha}  - \partial_{\beta} &  \partial_{\alpha}    
\end{array} \right] \; ; 
\end{equation}

\vspace{1\baselineskip}

\begin{equation}\label{hgam2}
\begin{array}{l}
\gamma^{\alpha} = \left[ \begin{array}{cccc}
0 & 0 & 1 & 0 \\
0 & 0 & 0 & 0 \\
h^2+r & h & 0 & 0 \\
h & 1 & 0 & 0  \end{array} \right] \quad , \quad 
\gamma^{\beta} = \left[ \begin{array}{cccc}
0 & 0 & 0 & 0 \\
0 & 0 & 1 & 0 \\
-h & 0 & 0 & 0 \\
-1 & 0 & 0 & 0  \end{array} \right] \quad , \\
\, \\
\gamma^{\gamma} = \left[ \begin{array}{cccc}
0 & 0 & 0 & 1 \\
0 & 0 & 0 & 0 \\
-h & -1 & 0 & 0 \\
0 & 0 & 0 & 0  \end{array} \right] \quad , \quad
\gamma^{\delta} = \left[ \begin{array}{cccc}
0 & 0 & 0 & 0 \\
0 & 0 & 0 & 1 \\
1 & 0 & 0 & 0 \\
0 & 0 & 0 & 0  \end{array} \right] \quad .
\end{array}
\end{equation}

\vspace{1\baselineskip}

\noindent
{\bf 2. \,}  ${\cal M}_h^{(2)}$  case:

\begin{equation}\label{hgam3}
(\hat{R}_h^{\epsilon})^{-1} = \left[ \begin{array}{cccc}
-h^2 & 0 & 0 & 1 \\
2h & -1 & 0 & 0 \\
2h & 0 & -1 & 0 \\
1 & 0 & 0 & 0  \end{array} \right] \quad ; \quad 
Y^{\epsilon} = \left[ \begin{array}{cc}
-h^2 \partial_{\alpha} + \partial_{\delta} & 2h \partial_{\alpha} 
- \partial_{\gamma}  \\
2h \partial_{\alpha}  - \partial_{\beta} &  \partial_{\alpha}    
\end{array} \right] \quad ;
\end{equation}

\vspace{1\baselineskip}

\begin{equation}\label{hgam4}
\begin{array}{l}
\gamma^{\alpha} = \left[ \begin{array}{cccc}
0 & 0 & 1 & 0 \\
0 & 0 & 0 & 0 \\
-h^2 & 2h & 0 & 0 \\
2h & 1 & 0 & 0  \end{array} \right] \quad , \quad 
\gamma^{\beta} = \left[ \begin{array}{cccc}
0 & 0 & 0 & 0 \\
0 & 0 & 1 & 0 \\
0 & 0 & 0 & 0 \\
-1 & 0 & 0 & 0  \end{array} \right] \quad , \\
\, \\
\gamma^{\gamma} = \left[ \begin{array}{cccc}
0 & 0 & 0 & 1 \\
0 & 0 & 0 & 0 \\
0 & -1 & 0 & 0 \\
0 & 0 & 0 & 0  \end{array} \right] \quad , \quad
\gamma^{\delta} = \left[ \begin{array}{cccc}
0 & 0 & 0 & 0 \\
0 & 0 & 0 & 1 \\
1 & 0 & 0 & 0 \\
0 & 0 & 0 & 0  \end{array} \right] \quad .
\end{array}
\end{equation}
\noindent
Let us now discuss the plane wave solutions of the $h$-K-G operator.

\setcounter{equation}{0}

\section{$h$-deformed invariant equations: solutions of the K-G equation}

\indent
   
   Relativistic invariant equations have been already discussed for the 
$q$-deformation \cite{SONG,SCHI,PILLIN,FTUV94-21,MeyerM}
and recently the case of the deformed Klein-Gordon and Dirac equations on
a deformed homogeneous space (Minkowski) has been considered in \cite{PODLES}.
We shall follow here a method similar to that followed in 
\cite{PLB95}. This means
that together with the coordinate and derivative algebras ${\cal M}_h^{(j)}$,
${\cal D}_h^{(j)}$ ($j$=1,2) we shall introduce another algebra,
the algebra of  {\it momenta}  ${\cal P}_h^{(j)}$. The generators of 
${\cal P}_h^{(j)}$  will be  
the elements  of a 2$\times$2 matrix $P$. We may look at $P$ as the eigenvalues
matrix of the derivative matrix (in a naive quantization we may set 
$P/i \hbar$$\sim$$Y$ where $\hbar$
is the Planck constant). Thus, $P$ transforms
as a `covariant' vector under an element of $L^{(j)}$,
$P \mapsto (M^{\dagger})^{-1} P M^{-1}$ (cf. (\ref{dif1})).

Since we make the natural assumption $P \propto Y$, 
it is equally natural to assume that the elements of $P$ neither commute among
themselves nor with those of $K$ and $Y$ (cf. \cite{PODLES}). The
requirement that these commutation relations are preserved under the
Lorentz coaction may be again expressed by reflection equations. As 
 mentioned in Sec. 2, their $R$-matrix structure in the
$h$-deformed case  is unique: their form depends only on the 
contravariant (eq. (\ref{2.6})), covariant  (eq. (\ref{dif2})) or mixed
(eq. (\ref{dif4})) character of the objects that they contain. Hence, 
the commutation relations of the elements in $P^{(j)}$ (generators of
${\cal P}_h^{(j)}$) with those of $K^{(j)}$ (id. of ${\cal M}_h^{(j)}$)
and $Y^{(j)}$ (id. of ${\cal D}_h^{(j)}$) are fully 
determined by the equations

\begin{equation}\label{p1}
R_h^{\dagger}P_2R^{(2) \,-1}P_1=P_1R^{(3) \,-1}P_2R_h \quad,
\end{equation}
\begin{equation}\label{p2}
P_2R_hK_1R^{(2)}=R^{(3)}K_1R^{\dagger}_hP_2  \quad,
\end{equation}
\begin{equation}\label{p3}
R_h^{\dagger}Y_2R^{(2) \,-1}P_1=P_1R^{(3) \,-1}Y_2R_h \quad.
\end{equation}
\noindent
As the commutation relations among the elements of $P$ themselves are the
same as those of the derivatives $Y$, it follows that we can take $K$ and $P$
simultaneously hermitian and that 
\begin{equation}\label{p4}
P^2 := \frac{1}{2+h^2} tr_h(P^{\epsilon}P) \quad, \quad P^{\epsilon} :=
(\hat{R}_h^{\epsilon})^{-1}P
\end{equation}
\noindent
is central in the ${\cal P}_h^{(j)}$ algebra (compare with the definition of 
the d'Alembertian, eq. (\ref{dif3})) and in ${\cal M}_h^{(j)}$, 
${\cal D}_h^{(j)}$.

The first step in the search for  solutions of the deformed equations 
is to see whether the standard plane wave solutions
 $e^{-ipx/\hbar}$ are  modified for the $h$-deformed
Klein-Gordon equation. To see that they are not, we first introduce 
the $h$-Lorentz
$L_h^{(j)}$  invariant scalar product for coordinates and momenta
\begin{equation}\label{p5}
(K,P):=tr_h(KP) \quad.
\end{equation}
\noindent
It is not difficult to check now, using the invariance of the $h$-trace 
(\ref{dilat3}) and eqs. (\ref{p1}) and (\ref{p2}), that $(K,P)$ commutes
with the elements of the ${\cal M}_h^{(j)}$ and ${\cal P}_h^{(j)}$
algebras,
\begin{equation}\label{p6}
[(K,P),K]=0=[(K,P),P]  \quad.
\end{equation}
\noindent
Similarly, and using that $tr_{h(1)}(R_{12}{\cal P})=I$ for the matrix $R_h$,
we find
\begin{equation}\label{p7}
Y (K,P) = (K,P)Y+P \quad.
\end{equation}
\noindent
Iterating this equation we obtain  for the $n$-th power of the scalar product
\begin{equation}\label{p8}
Y (K,P)^n = (K,P)^n Y  + n(K,P)^{n-1}  P \quad
\end{equation}
\noindent
so that, on a constant, eq. (\ref{p8}) gives the usual rule
\begin{equation}\label{p9}
Y (K,P)^n =  n(K,P)^{n-1}  P \quad.
\end{equation}
\noindent
Hence, using the {\it ordinary} expansion $\exp i(K,P)=\sum_{n=0}^{\infty}
\frac{1}{n!}(iK,P)^n$ we check that 
\begin{equation}\label{p10}
Y \, \exp i(K,P) = iP \, \exp i(K,P) \quad.
\end{equation}
\noindent
Taking again the derivative, now using $Y^{\epsilon}$, we finally
obtain
\begin{equation}\label{p11}
\Box_h \, \exp i(K,P) = - P^2 \, \exp i(K,P) \quad ,
\end{equation}
\noindent
where we recall that  $\Box_h$ [$P^2$] is given in (\ref{dif3}) [(\ref{p4})].
Thus, the deformed plane-wave solutions of the $h$-K-G equation have the same
{\it form} as  in the free case, and a mass-shell $P^2=m^2c^2$ may
be defined similarly.
A similar study could be performed for the Dirac equation since the 
expressions of the gamma matrices are known and the exponential 
$\exp i(K,P)$ commutes with $P$, so one can reduce the $h$-Dirac equation
to a matrix one with non-commutative entries $P_{ij}$, but we shall 
not do it here.

\setcounter{equation}{0}

\section{Conclusions} 

\indent

 We have studied  the $h$-deformations of the Lorentz 
group associated with $SL_h(2)$ and the corresponding Minkowski spacetime 
algebras. The twisted character of these deformations has been shown by 
exhibiting the twisted nature of $L_h^{(1,2)}$, and it is also measured
by  the $h$-dilatation operator $s$. Our analysis shows that the
irreducible representations of 
deformed Minkowski algebras are quite different for different deformations.

We have also studied the kernels of the $h$-deformed relativistic equations and
given the defining $h$-Clifford algebra relations and explicit form for
the $h$-Dirac matrices. The analysis of the solutions of the $h$-K-G
equation is facilitated by the fact that the $h$-Minkowski differential 
calculi, which are based on the twisted $SL_h(2)$ deformation, have a real
structure. 
We now conclude  by comparing our procedure and solutions with those 
in \cite{PODLES}. There, the deformed  momenta are also introduced as the
generators of an algebra (of exchange type for a triangular $R$-matrix).
However, and in contrast with our ${\cal P}_h^{(j)}$ algebra, these
momenta commute with the coordinates and the derivatives. This unnatural
simplification
of the commutation properties, which would correspond to replacing eqs.
(\ref{p2}), (\ref{p3}) by $P_2K_1=K_1P_2$ and $Y_2P_1=P_1Y_2$, produces
a more complicated action of the derivatives on the scalar product and
(at least within the present scheme) appears to be inconsistent 
with covariance, since the previous expressions lack  the 
$R$-matrices necessary to reorder the
non-commuting $M$ matrices which transform $P$, $K$ and $Y$.


\newpage

\noindent
{\bf Acknowledgements.} This research has been partially sponsored by the
Spanish  CICYT and DGICYT
(AEN 96-1669, PR 95-439). One of the authors (J.A.)
wishes to thank the hospitality extended to him at DAMTP and 
at St. John's College.  P.K. wishes to thank the Generalitat Valenciana
for supporting his stay in  Valencia University, where this work was started.

\appendix

\setcounter{equation}{0}

\section{Appendix: The 16$\times$16 $h$-Lorentz R-matrices}


\indent

The explicit expressions of the 16$\times$16 ${\cal R}$-matrices for 
$L_h^{(1,2)}$  are obtained from  eq. (\ref{ex5}),
where $R^{(1)}=R_h$ (eq.  (\ref{1.5})), $R^{(2)}={\cal P}R^{(3)}{\cal P}$,
the $R^{(3)}$ are given in (\ref{4.8}) and $R^{(4)}=R_h^{\dagger}$  
(cf. eq. (\ref{2.6})).
Thus,

\vspace{1\baselineskip}

\noindent
{\bf 1. \,}  ${\cal M}_h^{(1)}$  case:

\vspace{1\baselineskip}

\footnotesize

\noindent
{\normalsize ${\cal R}_h^{(1)} =$ } 
\begin{equation}\label{appe1}
\hspace{-2cm} \left[ 
\begin{array}{cccccccc}                                                 
     
     1  &   h &  - h & - h^2 & - h & h^2 &    h^2+r & h(r-h^2)  \\

     0  &   1  &   0  &   - h & 0 &    h  &   0   &  r-h^2  \\             

     0  &   0  &   1 &     h &   0  &   0  &   -h  &    h^2  \\

    0  &   0  &   0 &    1   &  0  &   0  &   0  &    h  \\

    0  &   0 &    0  &   0 &    1  &   -h   &  -h   &  h^2  \\

    0 &    0  &   0 &    0 &    0  &   1   &  0   &  -h  \\

    0  &   0  &   0  &   0 &    0 &    0 &    1  &   -h \\

    0   &  0 &    0 &    0 &    0  &   0  &   0  &   1  \\
 
      0  &   0   &  0  &   0 &    0   &  0  &   0  &   0 \\

    0 &    0    & 0   &  0   &  0  &   0  &   0   &  0  \\

    0 &    0 &     0   &  0 &    0  &   0  &   0  &   0 \\

    0  &   0  &   0  &   0 &    0   &  0  &   0   &  0  \\

    0  &   0  &   0  &   0  &   0  &   0    & 0   &  0  \\

    0 &    0  &   0  &   0   &  0  &   0    & 0    & 0 \\

    0 &    0  &   0  &   0   &  0   &  0    & 0   &  0 \\

    0  &   0  &   0  &   0    & 0 &    0  &   0 &    0 

\end{array}  \right.  \hspace{3cm}
\end{equation}
$$
\hspace{4cm} \left. \begin{array}{cccccccc}     
                                                  
      h  & h^2-r & h^2 & h(h^2-r) & -h^2 & h(h^2-r) & h(r-h^2) & h^4-r^2 \\

         0 &    h &    0  &   h^2 &    0 &    h^2  &   0   &  h(h^2+r) \\
                                                      
         0  &   0 &    -h  & -h^2-r & 0  &   0   &  h^2   & -h(h^2+r)  \\
                                                       
         0  &   0   &  0  &   -h & 0  &   0  &   0  &   -h^2 \\
                                                    
         0  &   0 &    0  &   0  &    h  & -h^2-r & h^2 & -h(h^2+r) \\

         0  &   0  &   0  &   0   &  0 &   h &    0 &    h^2 \\
                                                       
         0   &  0 &    0  &   0   &  0  &   0 &  -h  &  h^2-r \\

       0 &    0 &    0 &    0   &  0  &   0  &   0   &  -h  \\

      1  &   h  &   h &    h^2   &  -h & h^2   & r-h^2 & h(h^2+r) \\
                                                      
         0  &   1  &   0   &  h    & 0  &   h  &   0  &  h^2+r  \\

         0  &   0  &   1   &  h   &  0  &   0   &  -h   &  h^2 \\

         0   &  0  &   0  &   1 &    0 &    0   &  0  &    h \\
                                                       
         0  &   0 &    0  &   0  &   1  &   -h &  h  & -h^2 \\

         0   &  0 &    0 &    0  &   0  &   1 &    0  &    h \\

         0   &  0  &   0 &    0   &  0  &   0  &   1    & -h  \\

         0  &   0  &   0  &   0   &  0   &  0  &   0    & 1
\end{array}
 \right]
$$

\newpage

\normalsize

\noindent
{\bf 2. \,}  ${\cal M}_h^{(2)}$  case:

\vspace{1\baselineskip}

\small

\noindent
{\normalsize ${\cal R}_h^{(2)} =$ }    
\begin{equation}\label{appe2}
\left[ \begin{array}{cc}
\begin{array}{cccccccc}                                                 
     
     1  &   2 h &  -2 h & -4 h^2 & -2 h & 0 &    -h^2 &   -2 h^3  \\

     0  &   1  &   0  &   -2 h & 0 &    0  &   0   &  -h^2  \\             

     0  &   0  &   1 &    2 h &   0  &   0  &   0  &   4 h^2  \\

    0  &   0  &   0 &    1   &  0  &   0  &   0  &   2 h  \\

    0  &   0 &    0  &   0 &    1  &   0   &  0   &  0  \\

    0 &    0  &   0 &    0 &    0  &   1   &  0   &  0  \\

    0  &   0  &   0  &   0 &    0 &    0 &    1  &   0 \\

    0   &  0 &    0 &    0 &    0  &   0  &   0  &   1 \\

    0  &   0   &  0  &   0 &    0   &  0  &   0  &   0 \\

    0 &    0    & 0   &  0   &  0  &   0  &   0   &  0  \\

    0 &    0 &     0   &  0 &    0  &   0  &   0  &   0 \\

    0  &   0  &   0  &   0 &    0   &  0  &   0   &  0  \\

    0  &   0  &   0  &   0  &   0  &   0    & 0   &  0  \\

    0 &    0  &   0  &   0   &  0  &   0    & 0    & 0 \\

    0 &    0  &   0  &   0   &  0   &  0    & 0   &  0 \\

    0  &   0  &   0  &   0    & 0 &    0  &   0 &    0 
  
\end{array} & \hspace{-0.5cm} \begin{array}{cccccccc}     
                                                  
         2 h  & h^2   &  0  &   0   &  -4 h^2 & 0  &   -2 h^3 & 3 h^4 \\

         0 &    0 &    0  &   0 &    0 &    0  &   0   &  0 \\
                                                      
         0  &   0 &    0  &   -3 h^2 & 0  &   0   &  0   &  -6 h^3  \\
                                                       
         0  &   0   &  0  &   -2 h & 0  &   0  &   0  &   -4 h^2 \\
                                                    
         0  &   0 &    0  &   0  &   2 h  & -3 h^2 & 4 h^2 &  -6 h^3 \\

         0  &   0  &   0  &   0   &  0 &    0 &    0 &    0 \\
                                                       
         0   &  0 &    0  &   0   &  0  &   0 &    0  &   -3 h^2 \\

       0 &    0 &    0 &    0   &  0  &   0  &   0   &  -2 h  \\     
                            
         1  &   0  &   0 &    0   &  -2 h & 0   &  -h^2 &   0 \\
                                                      
         0  &   1  &   0   &  0    & 0  &   0  &   0  &   3 h^2  \\

         0  &   0  &   1   &  0   &  0  &   0   &  0   &  0 \\

         0   &  0  &   0  &   1 &    0 &    0   &  0  &   2 h \\
                                                       
         0  &   0 &    0  &   0  &   1  &   -2 h & 2 h  & -4 h^2 \\

         0   &  0 &    0 &    0  &   0  &   1 &    0  &   2 h \\

         0   &  0  &   0 &    0   &  0  &   0  &   1    & -2 h  \\

         0  &   0  &   0  &   0   &  0   &  0  &   0    & 1
\end{array}
\end{array}  \right]
\end{equation}

\end{document}